# ON FINITE TYPE 3-MANIFOLD INVARIANTS IV: COMPARISON OF DEFINITIONS

STAVROS GAROUFALIDIS AND JEROME LEVINE

This edition: September 12, 1995; First edition: September 12, 1995

Fax number for J. Levine: (617) 736 3085

Email: stavros@math.mit.edu

Email: levine@max.math.brandeis.edu

ABSTRACT. The present paper is a continuation of [Ga], [GL] and [GO]. Using a key Lemma we compare the two currently existing definitions of finite type invariants of oriented integral homology spheres and show that type $3m$ invariants in the sense of Ohtsuki [Oh] are included in type $m$ invariants in the sense of the first author [Ga]. This partially answers question 1 of [Ga]. We show that type $3m$ invariants of integral homology spheres in the sense of Ohtsuki map to type $2m$ invariants of knots in $S^3$, thus answering question 2 from [Ga].

## Contents

The authors were partially supported by NSF grants DMS-95-05106 and DMS-93-03489.

# 1. Introduction

## 1.1. Definitions.
We begin by recalling some definitions from [Ga], [GO] and [GL] and establishing some notation that will be followed in the present paper.

All 3-manifolds considered are oriented integral homology spheres. A link $L$ in a integral homology sphere is called *algebraically split* (denoted $AS$) if the linking numbers between its components vanish. A link $L$ is called boundary if each component bounds a Seifert surface, and the Seifert surfaces are disjoint from each other. A (integral) *framing* $f$ for a link $L$ in an integral homology sphere $M$ is a sequence of integers indicating the linking numbers of longitudes of $L$ with the corresponding components. This requires a choice of orientation, but if one gives the longitudes the parallel orientation then the framing number is independent of the choice of orientation. A link is called *unit-framed* if the framing on each component is $\pm 1$. A framed link $(L, f)$ is called $AS$-admissible (respectively, $B$-admissible) if it is $AS$ (resp. boundary) and unit-framed. Let $\mathcal{M}$ denote the $\mathbb{Q}$-vector space generated by the diffeomorphism classes of oriented integral homology 3-spheres. Let

$$[M, L, f] = \sum_{L' \subseteq L} (-1)^{|L'|} M_{L', f'} \tag{1}$$

where $f$ denotes a framing of $L$, $f'$ is the restriction of $f$ to $L'$ and $M_{L,f}$ denotes Dehn surgery on the framed (unoriented) link $L$ in $M$. $|L|$ denotes the number of components. Let $\mathcal{F}_n^O \mathcal{M}$ (resp. $\mathcal{F}_n^G \mathcal{M}$) be the subspace of $\mathcal{M}$ spanned by all $[M, L, f]$ for $AS$-admissible (resp. $B$-admissible) links $L$ of $n$ components in integral homology spheres $M$. Obviously, $\mathcal{F}_*^O \mathcal{M}, \mathcal{F}_*^G \mathcal{M}$ are decreasing filtrations on $\mathcal{M}$.

$v : \mathcal{M} \to \mathbb{Q}$ is a *type $m$* invariant of integral homology spheres if $v(\mathcal{F}_{m+1}^O \mathcal{M}) = 0$, see [Oh]. Similarly, $v$ is of $B$-type $m$ if $v(\mathcal{F}_{m+1}^G \mathcal{M}) = 0$, see [Ga]. We denote the space of type $m$ invariants of integral homology spheres by $\mathcal{F}_m \mathcal{O}$.

## 1.2. Statement of the results.

**Theorem 1.** *With the above notation we have:*

$$\mathcal{F}_n^G \mathcal{M} \subseteq \mathcal{F}_{3n}^O \mathcal{M} \tag{2}$$

**Corollary 1.1.** *Type $3m$ invariants of integral homology spheres are included in $B$-type $m$ invariants of integral homology spheres.*

**Corollary 1.2.** *If $\lambda \in \mathcal{F}_{3m} \mathcal{O}$, $K$ a knot in a integral homology sphere $M$, $n \in \mathbb{Z}$, then $\lambda(M_{K,1/n})$ is a polynomial of $n$ of degree $m$.*

**Theorem 2.** 
- *We have the following equality of filtrations:*

$$\mathcal{F}_{3n}^O \mathcal{M} = \cap_{k \geq 0} (\mathcal{F}_n^G \mathcal{M} + \mathcal{F}_k^O \mathcal{M}) \tag{3}$$

- *Assuming that for every $n \geq 0$ there is a $k \geq 0$ such that $\mathcal{F}_k^O \mathcal{M} \subseteq \mathcal{F}_n^G \mathcal{M}$, we obtain that $\mathcal{F}_{3n}^O \mathcal{M} = \mathcal{F}_n^G \mathcal{M}$.*

**Conjecture 1.** *For every $n \geq 0$ there is a $k \geq 0$ such that $\mathcal{F}_k^O \mathcal{M} \subseteq \mathcal{F}_n^G \mathcal{M}$.*

Recall from [Ga] that there is a well-defined map $\Phi : \mathcal{F}_n \mathcal{O} \to \mathcal{F}_{n-1} \mathcal{V}$ where $\mathcal{F}_n \mathcal{O}$ denotes the space of type $n$ invariants of integral homology spheres and $\mathcal{F}_{n-1} \mathcal{V}$ denotes the space of type $n - 1$ invariants of knots in $S^3$. In [Ga], it was asked whether $\Phi$ descends to a map $\mathcal{F}_{3n} \mathcal{O} \to \mathcal{F}_{2n} \mathcal{V}$. In [GrLi] it was shown that $\Phi$ descends to a map $\mathcal{F}_n \mathcal{O} \to \mathcal{F}_{n-2} \mathcal{V}$ if $n \geq 4$. In [GL] we showed that $\Phi$ descends to a map $\mathcal{F}_{5n+1} \mathcal{O} \to \mathcal{F}_{4n} \mathcal{V}$. Recently N. Habegger [Ha] gave a proof of the above question. We will give a different proof along the lines of our argument in [GL]. We first show the following theorem:

**Theorem 3.** *If $L$ is an AS-admissible link containing a $2m + 1$-component trivial sublink, then $[S^3, L, f] \in \mathcal{F}_{3m}^O \mathcal{M}$.*

As in [GL], Theorem 3 implies the following result:

**Theorem 4.** [Ha] *The above map $\Phi$ factors through a map:*

$$\tag{4} \mathcal{F}_{3m} \mathcal{O} \to \mathcal{F}_{2m} \mathcal{V}$$

### 1.3. Questions.

**Question 1.** Can every integral homology sphere be obtained by Dehn surgery on a unit-framed *boundary* link in $S^3$?

*Remark* 1.3. It is recently shown by [Au] and [GoLu] that there are integral homology spheres that cannot be obtained by surgery on a knot.

### 1.4. Plan of the proof.
In section 2 we prove a key Lemma 2.1. In section 3.1 we give a proof of Theorem 1 and corollaries 1.1 and 1.2. In section 3.2 we give a proof of Theorem 2. In section 4 we give a proof of Theorems 3 and 4.

### 1.5. Acknowledgement.
We wish to thank D. Bar-Natan for for many useful conversations. Especially we wish to thank the `Internet` for providing a continuous channel for conversations, information and occasional frustrations.

## 2. A key Lemma

This section is devoted to the proof of the following Lemma, which is the key to the proof of Theorems 1 and 3. Note that *all links considered in the rest of the paper are unit-framed.*

**Lemma 2.1.** *Let $L$ be an AS-admissible link containing a sublink with two components $k_1, k_2$ which bound disks $D_1, D_2$ in $U$ so that $D_1 \cap D_2$ is a single arc $\alpha$ in the interior of $D_1$ (a ribbon intersection). Suppose $U$ is a ball containing $D_1 \cap D_2$ whose* *intersection with $L$ is as pictured in Figure 1. Let $L_\alpha$ be the link obtained from $K$ by replacing $k_1$ with $k'_1$, a small circle in $D_1$ about $\alpha$. See Figure 2. Then:*

(5) $\qquad [S^3, L, f] = [S^3, L_\alpha, f] + \text{ a linear combination of } [S^3, L(\nu), f(\nu)]$

*where each link $L(\nu)$ contains $L$ as a proper sublink such that $L(\nu) - L \subseteq U$ (we will say such links are subordinate to $L$).*

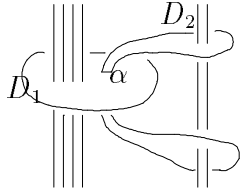

**Figure 1.** Shown here is the intersection of $L$ with $U$. Note that $k_i = \partial D_i$ for $i = 1, 2$. and that the discs $D_1$ and $D_2$ intersect in an ribbon arc $\alpha$.

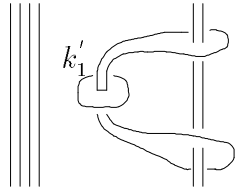

**Figure 2.** Shown here is the intersection of $U$ with the link $L_\alpha$ obtained by changing $k_1$ to $k'_1$. Note that $L_\alpha - L \subseteq U$.

*Proof.* Let $L_{twist}$ be the link obtained from $L$ by replacing $L \cap U$ with Figure 3. We first show:

**Claim 2.2.**

(6)
$[S^3, L_{twist}, f] = -[S^3, L, f] + 2[S^3, L_\alpha, f] + \text{ a linear combination of } [S^3, L(\nu), f(\nu)]$

*where the $\{L(\nu)\}$ are subordinate to $L$.*

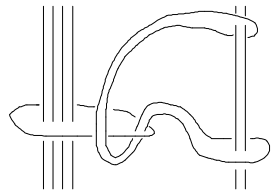

**Figure 3.** Shown here is the intersection of U with the link $L_{twist}$. Note that $L_{twist} - L \subseteq U$.

*Proof.* [of claim 2.2] Consider $L_{untwist}$ obtained from $L$ by replacing $L \cap U$ by Figure 4.

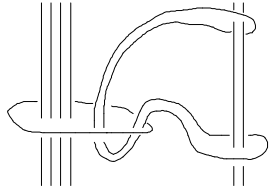

**Figure 4.** Shown here is the intersection of U with the link $L_{untwist}$. Note that $L_{untwist} - L \subseteq U$.

Now apply [GL][ Theorem 5] to the disc $D_1$, where we use three bands. The last two are the ones seen penetrating $D_1$ in Figure 4 and the first one contains all the other strands of $L_{untwist}$ penetrating $D_1$. Then Theorem 5 implies that $[S^3, L_{untwist}, f]$ is a sum of six terms in which $D_1$ is replaced by smaller subdisks and others in which $D_1$ is replaced by more than one subdisk. These last terms are all subordinate to $L$. Three of the first six terms are just $[S^3, L, f], [S^3, L_{twist}, f]$ and $[S^3, L_{notwist}, f]$, where $L_{notwist}$ is obtained from $L_{untwist}$ by replacing $D_1$ by a subdisk which only encloses the two penetrations of the third band $\beta$. But $[S^3, L_{untwist}, f] = [S^3, L_{notwist}, f] = 0$ because we can obviously isotop $\beta$ to miss $D_1$ and then $k_2$ bounds a disk in the complement of the rest of the link. Two of the remaining three terms are $-[S^3, L_\alpha, f]$ and the last term is given by a link obtained from $L$ by replacing $D_1$ by a subdisk disjoint from $\beta$ but intersected by all the other strands of $L$ which intersect $D_1$. As above this term vanishes, since $k_2$ bounds a disk in the complement of the rest of the link. □

**Claim 2.3.**

$$[S^3, L_{twist}, f] = [S^3, L, f] + \text{ a linear combination of } [S^3, L(\nu), f(\nu)]$$

where the $\{L(\nu)\}$ are all subordinate to $L$.

*Proof.* [of claim 2.3] After an isotopy, $L_{twist} \cap U$ appears as in Figure 5.

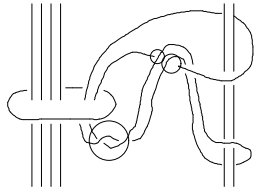

**Figure 5.** Shown here is an isotopy fixing the boundary of the intersection of U with the link $L_{twist}$. Note that $L_{twist} - L \subseteq U$. Also circled are 3 crossings to be changed.

Three crossing changes, from Figure 5, will convert this into $L \cap U$. These crossing changes are effected by surgeries along three circles. In Figure 6 we see $L \cap U$ with

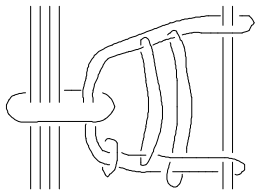

**Figure 6.** Yet another intersection of $U$ with a link.

the three circles added. Thus we conclude that:

$$(7) \quad [S^3, L_{twist}, f] = [S^3, L, f] + \text{ a linear combination of } [S^3, L(\nu), f(\nu)]$$

where the $L(\nu)$ consist of $L$ together with one or more of the extra circles in Figure 6. □

Obviously Lemma 2.1 follows from Claims 2.2 and 2.3. □

## 3. Proof of Theorems 1 and 2

This section is devoted to the proof of theorem 1 and 2.

**3.1. Proof of Theorem 1.** For the convenience of the reader, we divide the proof of Theorem 1 in 6 steps. We begin with some definitions that will be useful. A pair of links $(L, L_b)$ is called *n-boundary* if $L_b$ is a sublink of $L \subseteq S^3$ and $L_b$ is a boundary $n$ component link in the complement of $L - L_b$. The *goodness* $k(L, L_b)$ of an $n$-boundary pair is the number of components of $L - L_b$. The *genus* $g(L, L_b)$ of an $n$-boundary pair $(L, L_b)$ is the minimal total genus of disjoint Seifert surfaces of $L_b$ in the complement of $L - L_b$.

*Proof.* [of theorem 1] With the above terminology, we have the following step 1:

- **Step 1** $\mathcal{F}_n^G \mathcal{M}$ is generated by all $[S^3, L, f]$ for all $n$-boundary pairs $(L, L_b)$.

*Proof.* Let $\tilde{\mathcal{F}}_n^G \mathcal{M}$ denote the subspace spanned by all $[S^3, L, f]$ for all $n$-boundary pairs $(L, L_b)$. We first show $\mathcal{F}_n^G \mathcal{M} \subseteq \tilde{\mathcal{F}}_n^G \mathcal{M}$. Write $M = S^3_{L', \delta}$ for an algebraically split unit-framed link $L'$ in $S^3$. Since $L$ is an $n$-component boundary link in $M$, we can assume that $L$ bounds Seifert surfaces $\Sigma$ such that $\Sigma \cap L'$ is empty (here we mean by $L'$ the corresponding tubes of $M$). Thus $L \cup L'$ becomes a link in $S^3$, and $(L \cup L', L)$ is a $n$-boundary pair. We now proceed by upward induction on the number of components $|L'|$ of $L'$. If $L'$ is empty, we are done by definition. Otherwise, using Equation (1) we get:

$$(8) \quad [S^3, L \cup L', f \cup \delta] = \pm[M, L, f] + \sum_{L'' \subsetneq L'} \pm[S^3_{L'', \delta''}, L, f]$$

By induction, all the terms in the summation on the right hand side belong to $\tilde{\mathcal{F}}_n^G \mathcal{M}$ and so we conclude that $[M, L, f]$ does also.

The fact that $\tilde{\mathcal{F}}_n^G \mathcal{M} \subseteq \mathcal{F}_n^G \mathcal{M}$ is an immediate consequence of Equation (8). □

Let $(L, L_b)$ be an $n$-boundary pair. We want to show that $[S^3, L, f] \in \mathcal{F}^O_{3n}\mathcal{M}$. We proceed by *primary downward* induction on the goodness $k(L, L_b)$, and *secondary upward* induction on the genus $g(L, L_b)$. If $k(L, L_b) \geq 2n$ we are done by definition. If $g(L, L_b) = 0$ we are also done, since $[S^3, L, f] = 0$.

- **Step 2** We may assume that the components of $L - L_b$ are all unknotted.

*Proof.* This can be achieved by crossing changes in $L - L_b$ and, since is the result of a $\pm 1$-surgery along a small circle $C$ enclosing the crossing, the change to $[S^3, L, f]$ is given by an element $[S^3, L \cup C, f \cup \pm 1]$. See Figure 7. Since $(L \cup C, L_b)$ remains an $n$-boundary pair, whose goodness is one more than the goodness of $(L, L_b)$, it follows by the primary inductive hypothesis that $[S^3, L \cup C, f \cup \pm 1] \in \mathcal{F}^O_{3n}\mathcal{M}$. □

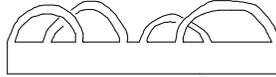

**Figure 7.** A relation in $\mathcal{M}$. Here an unknot circles the same component of a link, with linking number zero.

- **Step 3** Suppose that $L_b = \partial \Sigma_b$, where $\Sigma_b$ is a union of Seifert surfaces in the complement of $L - L_b$. We may assume that $\Sigma_b$ is embedded in a standard, almost planar (except for the necessary band crossings) way. See Figure 8.

**Figure 8.** A surface of genus 2 (and 4 bands) whose boundary is an unknot.

*Proof.* This can be achieved by band crossing changes, which are the result of a $\pm 1$-surgery along a circle enclosing the band crossing. Actually this surgery will introduce some extra twists into the bands, but further suregery along circles enclosing these twists will remove them. See Figure 9. As in step 2, the changes to $[S^3, L, f]$ are linear combinations of $[S^3, L', f']$ for $n$-boundary pairs $(L', L_b)$ with strictly higher goodness than that of $(L, L_b)$. By appealing to the primary inductive hypothesis, the changes to $[S^3, L, f]$ lie in $\mathcal{F}^O_{3n}\mathcal{M}$. □

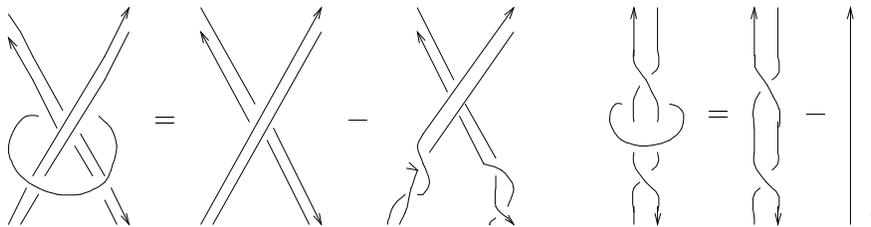

**Figure 9.** A few more identities in $\mathcal{M}$

Let $\{K_i\}$ denote the components of $L - L_b$. Since by step 2 they are unknotted, we may choose embedded disks $D_i$ so that $K_i = \partial D_i$. Furthermore, since $\Sigma_b$ is just a thickening of a wedge of circles, we may choose the $D_i$ so that their intersections with $\Sigma_b$ consist of a number of transverse penetrations of the interiors of the $D_i$ by the bands of $\Sigma_b$. See Figure 10. We will be interested in counting the number of "band penetrations".

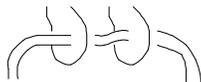

**Figure 10.** A band of a surface penetrating two pieces of discs.

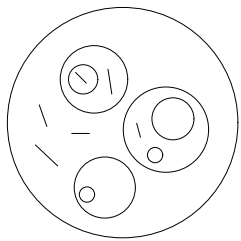

**Figure 11.** An intersection of a disc with the surface $\Sigma_b$

- **Step 4** We may assume that every band of $\Sigma_b$ penetrates at least one $D_i$.

*Proof.* Suppose that a band $\beta$ from one of the components $\Sigma$ of $\Sigma_b$ penetrates no $D_i$. We will show how to replace $\Sigma$ by a surface of lower genus and then appeal to the secondary inductive hypothesis. This will also involve a number of changes to L, but using the primary inductive hypothesis each of these changes will only be by an element of $\mathcal{F}_{3n}^O \mathcal{M}$.

Let $C$ be the circle in $\Sigma$ which goes once around the band $\beta$. $C$ bounds an obvious disk in the plane containing $\Sigma$. We push this disk slightly off the plane (except on $C$) to obtain a disk $D$ such that $D \cap \Sigma = \emptyset$ and $\partial D = C$. Now, since $C \cap \bigcup_i D_i = \emptyset$, $D \cap \bigcup_i D_i$ consists of circles and interior arcs, see Figure 11. If $D \cap \bigcup_i D_i = \emptyset$, then we can perform a surgery on $\sigma$ along $D$ to obtain the desired surface of lower genus. Thus we only have to see how to remove these intersections. We claim that we can first remove the arc intersections and then (by using an innermost circle argument) remove the circle intersections. In fact we only have to remove the arcs since since it is really only necessary that $D \cap (L - L_b) = \emptyset$. Suppose $\alpha$ is an arc and a component of $D_j \cap D$. We can perform an isotopy of $D_j$ to move $\alpha$ adjacent to the boundary $C$ of $D$. This may require $\alpha$ to cross some circle components of $D \cap \bigcup D_i$ which means that $D_j$ may cut through some $D_i$ during the isotopy. If $i = j$ the result will be that $D_j$ is now only immersed, but this will not be important. We have only a *regular homotopy* of $D_j$ but still an isotopy of $K_j$. Now a neighborhood of $\alpha$ in $D_j$ is a band which is adjacent to the band $\beta$. If we change this band crossing, the result will be to eliminate $\alpha$. As above this crossing change can be produced by a $\pm 1$-surgery on a small circle enclosing the two bands and so the change in $[S^3, L, f]$ is, by primary induction, an element of $\mathcal{F}^O_{3n}\mathcal{M}$. Note that we have not changed $\Sigma$. See Figure 12. □

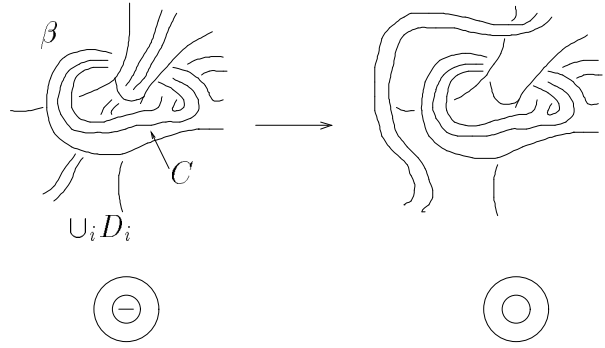

**Figure 12.** In the left side of this picture is shown a band $\beta$ that does not penetrate any of the discs $D_i$, its associated circle $C$. The disc which the circle $C$ bounds intersects the union $\cup_i D_i$ as shown in the lower part of the left hand side. After a band change move, shown on the right hand side, we can arrange so that the new disc of the new circle has one less band intersection with the union of the discs $\cup_i D_i$.

- **Step 5** We may assume that each disc $D_i$ has at most two band penetrations.

*Proof.* We want to apply [GL, Theorem 5] to every disc $D_i$. The bands in Theorem 5 are all but one the bands of $\Sigma_b$ penetrating $D_i$ and the remaining band consists of all the strands of $L - L_b$ penetrating $D_i$. Thus $[S^3, L, f]$ is a sum of elements in which $D_i$ is replaced by one or more disks inside $D_i$ containing no more than two bands of $\Sigma_b$. Thus $L - L_b$ is changed, but not $\Sigma_b$ and so, the change in $[S^3, L, f]$ comes from $n$-boundary pairs of higher goodness than $(L, L_b)$, and thus, by the primary inductive hypothesis, lie in $\mathcal{F}^O_{3n}\mathcal{M}$. □

Let $\{\Sigma_j | 1 \leq j \leq n\}$ denote the connected components of $\Sigma_b$.
- **Step 6** We may assume that, if $\Sigma_j$ has genus one and a band of $\Sigma_j$ penetrates only one disk $D_i$, then $D_i$ is penetrated by no other bands of $\Sigma_b$.

*Proof.* Suppose one of the bands $\beta$ of $\Sigma_j$ penetrates $D_i$ once. Let $L_\beta$ be defined from $L$ by replacing $k_i$ with a small meridian circle about $\beta$. Then Step 6 will be confirmed by the following:

**Claim 3.1.** *We have the following equality:*

(9) $$[S^3, L, f] = [S^3, L_\beta, f] \mod \mathcal{F}^O_{3n}\mathcal{M}$$

*Proof.* [of Claim 3.1] This Claim follows from the key Lemma 2.1 as follows. We can draw $L \cap U$, where $U$ is a ball containing $\Sigma_j$, as in Figure 13. If we expunge $\Sigma_j$ from the picture we have exactly the situation in Figures 1 and 2 of Lemma 2.1. If we put $\Sigma_j$ back into any of the $L(\nu)$ of Lemma 2.1, we see that it may intersect the additional components of $L(\nu)$. However we can add tubes to $\Sigma_j$ to eliminate these intersections and, since none of other $\Sigma_i$ intersect $U$, we may conclude from the key Lemma 2.1 that:

(10) $\qquad [S^3, L, f] = [S^3, L_\beta, f] +$ a linear combination of $[S^3, L(\nu), f(\nu)]$

where each of the pairs $(L(\nu), L_b)$ are $n$-boundary with strictly higher goodness than that of $(L, L_b)$. By the primary inductive hypothesis, we conculde the proof of Claim 3.1 and of step 5. □

□

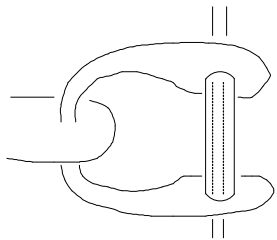

**Figure 13.** An intersection of $L$ with $U$. Shown also is the genus 1 surface $\Sigma_j$.

We can now complete the proof of Theorem 1. We define:

$r$ = number of $D_i$ penetrated by two bands
$s$ = number of $D_i$ penetrated by one band
$l$ = number of $\Sigma_j$ of genus $> 1$
$m$ = number of $\Sigma_j$ of genus $= 1$
$p$ = number of bands of surfaces of genus one penetrating only one disk

Obviously $k \geq r + s$ and $n = l + m$. By Step 4, the total number of band penetrations is $2r + s$. From Step 3 we thus conclude $2r + s \geq 4l + 4m - p$ and, by Step 5, we have $s \geq p$. Adding these two equations together gives us $2r + 2s \geq 4l + 4m$ and so:

$$k \geq r + s \geq 2l + 2m = 2n$$

This concludes the proof of Theorem 1. □

*Proof.* [of Corollary 1.1] If $v : \mathcal{M} \to \mathbb{Q}$ is of type $3m$, then $v(\mathcal{F}^O_{3m+1}\mathcal{M}) = 0$, and therefore, by Theorem 1, $v(\mathcal{F}^G_{m+1}\mathcal{M}) = 0$. □

*Proof.* [of Corollary 1.2] It follows by Exercise 4.2 of [Ga], using the remark that the $j^{th}$ cable (with zero framing) of a knot $K$ in a integral homology sphere $M$ is a boundary link of $j$ components. □

**3.2. Proof of Theorem 2.** In this section we prove theorem 2. The proof will use the $AS$ and $IHX$ relations on $\mathcal{M}$ proven in [GO]. For the convenience of the reader, we recall the notation and terminology from [GO]. A *Chinese manifold character* is a trivalent graph with vertex orientation. The degree of a Chinese manifold character is the number of edges of it. Let $\mathcal{CM}$ denote the vector space on the set of *Chinese manifold characters*, and let $\mathcal{BM}$ be the quotient space $\mathcal{CM}/\{AS, IHX\}$, where we quotient by the $AS$ and the $IHX$ relations of [GO]. We recall the following theorem from [GO]:

**Theorem 5.** [GO] *There is an onto map $O_m^\star : \mathcal{G}_m\mathcal{BM} \to \mathcal{G}_{3m}^O\mathcal{M}$.*

We need the following Lemma:

**Lemma 3.2.** *Let $\Gamma \in \mathcal{G}_{3m}\mathcal{BM}$ be a Chinese manifold character of $3m$ edges. Then $O_{3m}^\star(\Gamma) \in \mathcal{F}_{3m}^O\mathcal{M}$ actually lies in $\mathcal{F}_m^G\mathcal{M}$.*

*Proof.* Choose a circuit in $\Gamma$ (that is a sequence of edges, the beginning of which is the end of the previous, such that the end of the last is the beginning of the first, and such that the edges in the sequence are distinct). Color the edges of the circuit red. Thinking of the red colored edges of $\Gamma$ as the external circle, and using repeatedly the $IHX$ relation of [GO], we can write $O_m^\star(\Gamma)$ as a linear combination of values (under $O_m^\star$) of chord diagrams based on the red circle. By counting degrees, we see that each of the above mentioned chord diagrams have $m$ chords. Now using lemma 3.4 from [Ga] we see that the pairs $(O_m^\star(\text{chord diagram}), O_m^\star(m\text{-chords}))$ is a $m$-boundary pair, from which our conclusion follows. □

*Proof.* [of theorem 2] We can now finish the proof of theorem 2 as follows: Theorem 5 and Lemma 3.2 show that $\mathcal{F}_{3n}^O\mathcal{M} = \mathcal{F}_n^G\mathcal{M} + \mathcal{F}_{3n+1}^O\mathcal{M}$. Iterating the above equation we obtain Equation 3. The proof of theorem 2 is complete. □

## 4. Proof of Theorem 3

*Proof.* [of Theorem 3] Let $(L, f)$ be an $AS$-admissible link containing a trivial sublink $L_{trivial}$ of $2m + 1$ components. We will use *downward* induction on the number $r$ of components of $L - L_{trivial}$ to show first that $[S^3, L, f] \in \mathcal{F}_{3m+2}^O\mathcal{M}$. Obviously, if $r \geq m + 1$, we are done. We refer the reader to the proof of Theorem 7 in [GL] for the first part of the argument. Let us denote the components of $L_{trivial}$ by $\{L_i\}_{i=1}^{2m}$. Then $L_i = \partial D_i$, where the $\{D_i\}$ are disjoint disks. We showed in [GL] that we may assume that $L - L_{trivial}$ consists of components $\{l_k\}$ such that each $l_k$ is either of the form $\sigma_{ij}, i \neq j$, (where $\sigma_{ij}$ is pictured in Figure 15 of [GL]), or a band sum of two $\sigma_{ij}$. We will refer to $l_k$ as *simple* in the former case and *composite* in the latter case.

12                STAVROS GAROUFALIDIS AND JEROME LEVINE

Note that $\sigma_{ij}$ intersects the disks $D_i$ and $D_j$, but no others. We will say that $L_i$ is *k-special* if the only the only component $l_s$ of $L_{trivial}$ intersecting $D_i$ is $l_k$ , and if $l_k = \sigma_{ij} \sharp \sigma_{rs}$ then $i \neq r, s$.

We now use Lemma 2.1 to make an important observation.

**Claim 4.1.** *We may assume that if $L_i$ is k-special then $l_k$ is simple.*

*Proof.* [of Claim 4.1] Suppose that $l_k$ is composite. Then there is a ball $U$ which intersects $L$ as in Figure 13. But we can redraw this so that it looks like Figure 1 of Lemma 2.1, with the two component $(k_1, k_2)$ distinguished sublink of $L$ being $(l_k, L_i)$. For the subordinate links $L(\nu)$ of Lemma 2.1, we see that $[S^3, L(\nu), f(\nu)] \in \mathcal{F}^O_{3m}\mathcal{M}$ by induction. Thus, using Lemma 2.1, we can assume that each $l_k$ is simple. □

We now complete the proof of Theorem 3 by a counting argument. We define:

$$\begin{aligned}
a &= \text{number of simple } l_k \\
b &= \text{number of composite } l_k \\
c &= \text{number of } k_i \text{ which are } k\text{-special for some } k \\
2m &= |L_{trivial}| \\
d &= 2m + 1 - c \\
r &= |L - L_{trivial}|
\end{aligned}$$

Obviously $r = a + b$. As pointed out in [GL] , we may as well assume that every $D_i$ is intersected by at least one $l_k$ (or else $[S^3, L, f] = 0$). Counting intersections of the $\{l_k\}$ with the $\{D_i\}$, we have $2a + 4b \geq c + 2d$. From Claim 4.1 we obtain the inequality $2a \geq c$. Adding these last two inequalities we get $4a + 4b \geq 2c + 2d$ or $2r \geq 2m + 1$ or $r \geq m + 1$.

This concludes the proof that $[S^3, L, f] \in \mathcal{F}^O_{3m+2}\mathcal{M}$. Using Corollary 3.5 of [GL] (see also Corollary 1.6 of [GO]) we deduce that $\mathcal{F}^O_{3m+2}\mathcal{M} = \mathcal{F}^O_{3m}\mathcal{M}$, which concludes the proof of Theorem 3. □

*Proof.* [of Theorem 4] It follows *verbatim* as in Proposition 3.9 of [GL], using Theorem 3 of the present paper. □

DEPARTMENT OF MATHEMATICS, MASSACHUSETTS INSTITUTE OF TECHNOLOGY, CAMBRIDGE, MA 02139
*E-mail address*: stavros@math.mit.edu

DEPARTMENT OF MATHEMATICS, BRANDEIS UNIVERSITY, WALTHAM, MA 02254-9110
*E-mail address*: levine@max.math.brandeis.edu